\documentclass[12pt]{article}
\usepackage{latexsym}
\usepackage{amssymb}
\newcommand{\OG}{\mathrm{O}}
\newcommand{\Fit}{\mathrm{Fit}}
\newcommand{\Cent}{\mathrm{C}}
\begin{document}
\newtheorem{lem}{Lemma}
\newtheorem{theo}{Theorem}
\newtheorem{defn}{Definition}
\newtheorem{prop}{Proposition}
\newtheorem{cor}{Corollary}
\newtheorem{ques}{Question}
\newcommand{\ind}{\mathrm{Ind}}
\bibliographystyle{alpha}
\title{The influence of conjugacy class sizes on the structure of finite groups: a survey}
\date{25 March 2010}
\author{A.R.Camina \& R.D.Camina \\ School of Mathematics, University of
East Anglia,\\ Norwich, NR4 7TJ, UK; a.camina@uea.ac.uk\\
Fitzwilliam College, Cambridge, CB3 0DG, UK;\\
rdc26@dpmms.cam.ac.uk}
\maketitle
\section{Introduction}\label{intro}

In this survey we consider the influence of the sizes of conjugacy
classes on finite groups. Over the last 30 years there have been
many papers on the topic and it would seem to be a good idea to try
to bring some of the key results together in one place. This is
especially relevant as some authors seem unaware of others writing
in the field as well as some of the older results which seem to get
reproved quite regularly. It is hoped that in writing this, less time
will be spent in reproving old results, enabling more progress to
be made on some of the more interesting problems.

How much information can one expect to obtain from the sizes of
conjugacy classes? Sylow in 1872 examined what happened  if there was
information about the sizes of all conjugacy classes, whereas in 1904
Burnside showed that strong results could be obtained if there was particular
information about the size of just one conjugacy class. Landau in 1903
bounded the order of the group in terms of the number of conjugacy classes
whilst in 1919 Miller gave a detailed analysis of groups with very few
conjugacy classes. Very little then seems to have been done until 1953 when both
Baer and It\^o published papers on this topic but with different conditions
on the sizes.

By looking at these early results it can been seen that much will depend on how much information is given
and it is important to be explicit. For example if one knows that
there is only one conjugacy class size then the group is abelian,
but this can be any abelian group. However if you know the
collection of conjugacy class sizes, that is the multiplicities,
then the order of the group is also known. However it would still
not be possible to identify the group. Some authors have considered
the situation where the multiplicities of the conjugacy class sizes
are used if the size is not $1$. This is particularly true when the
authors have been studying aspects of the problem related to graphs.
Again if we only demand information about the sizes of conjugacy
classes and not their multiplicities, the group $G$ and $G\times A$
will have the same set whenever $A$ is an abelian group.
So we can only state results modulo a direct abelian factor.

Another reason for examining conjugacy classes is their fundamental role in
understanding the group ring, especially over $\mathbb C$. Recall that if
$G$ is a group then the group ring ${\mathbb C}[G]$ is the vector space
over $\mathbb C$ generated with basis the elements of $G$. For each conjugacy
class $K$ of $G$ define the element $c_K=\sum_{g\in K} g$. Then $c_K$ is in
the centre of ${\mathbb C}[G]$ and the elements $c_K$ are a basis for the
centre. This gives rise to a complete link to the character theory of $G$.
If $K_1,K_2, \ldots, K_k$ are the conjugacy classes of $G$ then there are
integers $a_{rst}: 1\leq r,s,t, \leq k$ so that
$$c_{K_r}c_{K_s}=\sum_t a_{rst}c_{K_t}.$$
Knowledge of the $a_{rst}$  is equivalent to knowing the character\label{lam}
table, this is a result that goes back to Frobenius, see
Lam \cite[Theorem 7.6]{MR1606416}. So at this level of information we see
that knowing character and conjugacy class data are the same.

As we are not mainly concerned with characters, and only mention results
as comparative results we refer the reader to \cite{MR1645304} or
\cite{MR2270898}, amongst many others, for a discussion of character theory and some of the necessary definitions. We discuss these connections further in Section~\ref{char}.

We will often refer to {\em the index} of an element, this is just
the size of the conjugacy class containing the element. The benefit
of this definition is entirely linguistic. Given an element $g$ in
some group $G$ we can talk about the index of $g$ rather than
talking about the size of the conjugacy class containing $g$. So if
we are referring to elements we use the term index but if we are
talking about conjugacy classes we refer to size.

In Section~\ref{basics} we introduce the basic definitions and
consider some elementary results which are used time and time again,
and proved time and time again in many papers.

Section~\ref{AP} considers results which might loosely be described
as placing arithmetical conditions on the indices of elements. One
of the most famous is Burnside's $p^a$ theorem. In this the
restriction is that there is an index which is a prime-power. We
then look at some conditions which imply solubility, here the
results are quite weak. But we see a distinction between demanding
conditions on all the indices of elements and only on some. This
occurs quite frequently, so for example, if we know that the indices
of all elements are powers of a  given prime $p$ then the group is
essentially a $p$-group. But then we can ask what happens if we only
ask that $p$-elements have index a power of $p$. Here we are
demanding fewer restrictions on the indices but that we can recognise
the order of elements in the group.

The conjugate type vector is introduced in Section~\ref{CTP}. This
considers the indices in descending order and looks at various
properties of the group which can be deduced from this sequence.
Some authors have considered variations. For example for a fixed
prime $p$ they consider the sequence of indices for $p'$-elements.

There is a large body of work considering the influence of conjugacy
class sizes on $p$-groups. We do not mention these results here,
unless they fit in a natural way, but refer the reader to a survey
paper on $p$-groups by A. Mann; one section is devoted to
`Representations and Conjugacy Classes'~\cite{MR1716701}.

In Section~\ref{graphs} various graphs that can be constructed from the sets
of conjugacy class sizes are defined. The properties of the graphs and the
relation to the structure of the groups is examined. This has been a very
active area in recent years. This has also been the case where
character degrees have been used similarly, see the survey by Lewis
\cite{MR2397031}.

Examining the influence of the number of conjugacy classes is considered in
Section~\ref{noclasses} and, as mentioned previously, connections with
character theory are considered in Section 7.\\[2ex]

{\bf Notation:} The notation we use is standard. Let $G$ be a finite
group and $x$ an element of $G$. The centraliser of $X\subseteq G$ is
$\{ g \in G : xg = gx\,\forall x\in X \}$ and denoted by $\Cent_G(X)$, note if $X=\{x\}$ we drop the brackets. The conjugacy
class of $x$ in $G$ is denoted by $x^G$.  The index of $x$ will be denoted by $\ind_G(x)$.
We denote the derived group of $G$ by $G'$ and the Fitting subgroup of $G$ by
$\Fit(G)$. We comment, given a prime $p$, that the terms $p$-regular
and $p$-singular are often used. Note that an element being
$p$-regular is the same as saying that the element is a
$p'$-element, that is that its order is not divisible by $p$. A
$p$-singular element is one which is not a $p'$-element.  Also if $\pi$ is a set of primes then $\pi'$ is the
set of primes not in $\pi$. $\OG_{\pi}(G)$ is the largest normal subgroup of $G$ whose order is a $\pi$-number. Also $\OG^{\pi}(G)$
is the smallest normal subgroup whose factor group $G/\OG^{\pi}(G)$ is a
$\pi$-group. In a group $G$ a $p$-complement is a subgroup $H$ whose index
in $G$ is exactly the highest power of $p$ to divide the order of $G$.
We use CFSG to denote the Classification of Finite Simple Groups.
\section{Basic Definitions and Results}\label{basics}

Baer~\cite{MR0055340}  gave the following definition which we
will use:- \begin{defn}
 Let $G$ be a finite group and
let $x \in G$. The index of $x$ in $G$ is given by $[G:{\rm
C}_G(x)]$ and is denoted by $\mbox{Ind}_G(x)$.\end{defn} Note,
$\mbox{Ind}_G(x)$ is the size of the conjugacy class of $x$, by the
orbit-stabiliser theorem.
\begin{lem} Let $N$ be a normal subgroup of $G$. Then\\
(i) if $x \in N$, $\mbox{Ind}_N(x)$ divides $\mbox{Ind}_G(x)$.\\
(ii) if $x \in G$, $\mbox{Ind}_{G/N}(xN)$ divides $\mbox{Ind}_G(x)$.
\end{lem}

\begin{lem}\label{ecc1}\cite{MR0294481} If $p$ is a prime which does not divide $\ind_G(x)$
for all $x$ of $p'$-order then the Sylow $p$-subgroup of
$G$ is a direct factor of $G$.\end{lem} {\bf Proof.} Let $P$ be a Sylow
$p$-subgroup of $G$ and let $C=\Cent_G(P)$. We now show that $PC$
contains a conjugate of every element in $G$. The lemma then follows
by an old argument of Burnside \cite[\S 26]{MR0069818}.

If $g \in G$ we can write $g=xy$ where $[x,y]=1$ and $x$ is a
$p$-element and $y$ is a $p'$-element. By the hypothesis we have
that $\Cent_G(y)$ contains a conjugate of $P$. So, conjugating as necessary,
we can ensure
that $P\subset \Cent_G(y)$. There exists $h\in \Cent_G(y)$ so
that $x^h\in P$. Hence $(xy)^h=x^hy \in PC$ as required.$\Box$

\begin{cor}\label{ecc} If $p$ is a prime which does not divide $\ind_G(x)$
for all $x$ then the Sylow $p$-subgroup of $G$ is an
abelian direct factor of $G$.\end{cor}

We comment here that this shows a dichotomy with the related theory
for character degrees. For if $G$ is a finite group with a normal
abelian Sylow $p$-subgroup then $p$ does not divide any character
degree. This is a result of It{\^o} \cite{MR0044528}, the converse is
also true, though the proof was not completed until 1986 with an
application of CFSG \cite{MR866772}.

Another useful lemma is the following.
\begin{lem}\label{prod} Let $x$ and $y$ be two elements of $G$
such that $\Cent_G(x)\Cent_G(y)=G$. Then $(xy)^G=x^Gy^G$.
\end{lem}
{\bf Proof.} Consider $x^gy^h$ with $g,h \in G$. Clearly $x^gy^h$ is conjugate to $x^{gh^{-1}}y$. We can write
$gh^{-1} = ab$ where $a \in \Cent_G(x)$ and $b \in \Cent_G(y)$. Thus $x^{gh^{-1}}y = x^{ab}y$ which is conjugate
to $x^ay^{b^{-1}} =xy$ as required.$\Box$
\begin{cor}Let $x$ and $y$ be two elements of $G$ such that
Ind$_G(x)$ and Ind$_G(y)$ are coprime then
$(xy)^G=x^Gy^G$.\end{cor}

This idea was discussed by Tchounikhin in 1930 \cite{MR0000}. In
this paper he shows that if there are three indices which are pairwise coprime then
the group is not simple.

We now give some definitions which are useful and provide some
unifying themes to our survey.

\begin{defn} Let $G$ be a finite group and $g \in G$. Define\\
(i) $\sigma_G(g)$ to be the set of primes dividing ${\rm
Ind}_G(g)$.\\
(ii)  $\sigma^*(G) = {\rm max} \{ |\sigma_G(g)| : g \in G \}$ and \\
(iii) $\rho^*(G) = \bigcup_{g \in G} \sigma_G(g).$

\end{defn}

We note that $\rho^*(G)$ is just the set of primes that divide
$G/Z(G)$ by Corollary~\ref{ecc}, this set is also known as the set of
{\it eccentric} primes.

In 1953  It{\^o} introduced the notion of a conjugate type vector:
\begin{defn} \cite{MR0061597} The conjugate type vector of a group $G$ is the
vector $(n_1, n_2, \ldots, n_r, 1)$ where $n_1 > n_2 > \ldots n_r
>1$ are the distinct indices of elements of $G$. The conjugate rank of $G$,
${\rm crk}(G)$, is given by $r$.
\end{defn}

\begin{defn} (i) A group $G$ is called Frobenius if $G$ can be written
as a product of two groups of coprime order $K$ and $C$ where $K$ is normal
and $\Cent_G(x)\subseteq K$ for all $1 \neq x \in K$.
$K$ is called the kernel of $G$ and $C$ the complement.
For more structure see \cite[V.8]{MR0224703}.\\
(ii) A group $G$ is called quasi-Frobenius if $G/Z(G)$ is
Frobenius.
\end{defn}
Let $G$ be quasi-Frobenius and let the pre-image of the kernel and
the complement  be $K$ and $C$ respectively. Then, if both $K$ and $C$
are abelian, the non-trivial conjugacy class sizes of $G$ are
$m=|C/Z(G)|$ and $n=|K/Z(G)|$. Note $\gcd(m,n)=1$.

\section{Arithmetical Properties}\label{AP}

\subsection{Prime-power index}

Perhaps the earliest results are those of Sylow~\cite{MR1509796}
which says that a group, all of whose indices are a power of a given
prime, has a non-trivial centre and of Burnside~\cite{Bu} which says
that if one index is a power of a prime then the group is not
simple. These are two classic results which give information about
the group from some arithmetical properties of the set of indices.

In 1990 Kazarin proved the following extension
to Burnside's result:
\begin{theo}\cite{MR1076932}\label{Kaz} Let $G$ be a finite group and let $x$ be an element
of $G$ such that $\ind_G(x)=p^a$ for some prime $p$ and integer
$a$. Then $\langle x^G\rangle$ is a soluble subgroup of
$G$.\end{theo}

Whilst the proof of Burnside's result depends on ordinary character theory, Kazarin's result
depends on modular character theory.
Using Theorem~\ref{Kaz} Camina \& Camina~\cite{MR1633180} were able to prove
that any element of prime-power index is in the second Fitting
subgroup. Flavell~\cite{MR1935499} has shown that whether an
element is in the second Fitting subgroup is determined by the
behaviour of two generator subgroups. This led to an interesting
discussion of how the indices of elements in two generator
subgroups can determine the index of an element in the whole
group, see~\cite{MR2564854}.

Baer \cite{MR0055340}\label{baerp} characterized all finite groups such that
every element of prime-power order has prime-power index. He then
went on to raise the question of the characterization of those
groups whose $q$-elements, for just one prime $q$, have prime
power index.
Camina \& Camina~\cite{MR1633180} introduced the following idea
based on Baer:

\begin{defn} Let $G$ be a finite group and let $q$ be a prime such that
$q$ divides $|G|$. Then $G$ is a $q$-Baer group, or equivalently
has the $q$-Baer property, if every $q$-element of $G$ has prime
power index.\end{defn} They then proved the following theorem:
\begin{theo} \cite{MR1633180} Let $G$ be a $q$-Baer group for some prime
$q$. Then\\
{\rm (a)}  $G$ is $q$-soluble with $q$-length 1, and\\
{\rm (b)} there is a unique prime $p$ such that each $q$-element
has $p$-power index.\\
Further, let $Q$ be a Sylow $q$-subgroup of
$G$, then\\
{\rm (c)}  if $p=q$, $Q$ is a direct factor of $G$, or\\
{\rm (d)} if $p \neq q$, $Q$ is abelian, ${\rm{O}}_p(G)Q$ is
normal in $G$ and $G/{\rm{O}}_{q'}(G)$ is soluble.\end{theo}

  Berkovic \& Kazarin prove some very similar
results but also some different ones. One of their results is the
following:
\begin{theo}\cite{MR2111209} If the index of every $p$-element of order $p$ (if $p$
is odd) or 4 (if $p=2$) is a $p$-power, then $G$ has a normal $p$-complement.
\end{theo}

Beltr{\'a}n \& Felipe have also proved some similar results, \cite{MR2051367}.

Recall, Lemma~\ref{ecc1} says that if all $p'$-elements have
$p'$-index then the Sylow $p$-subgroup of $G$ is a direct factor.
An interesting strengthening of this result is due to Liu, Wang \& Wei:
\begin{theo} \cite{MR2111851} If $p$ is a prime which does not divide $\ind_G(x)$
for all $p'$-elements $x$ of prime-power order then the Sylow $p$-subgroup of
$G$ is a direct factor of $G$.
\end{theo}
The proof uses the result of \cite{MR636194} which shows that in a
transitive permutation group there is an element of prime-power
order which acts fixed-point-freely, this result uses
CFSG.\label{CFSG}

In another variation on the theme Dolfi and Lucido say a finite
group $G$ has property $P(p,q)$ if every $p'$-element has
$q'$-index. The inspiration for this definition came from ideas in
character theory. In particular, a group $G$ has the property
$BP(p,q)$ if every $p$-Brauer character has degree prime to $q$.
Dolfi \& Lucido prove (amongst other things) the
following:\label{DL}
\begin{theo} \cite{MR1826493} Let $G$ be a finite group satisfying $P(p,q)$ with
$p\neq q$. Then $\OG^p(G)$ is $q$-nilpotent and $G$ has abelian
Sylow $q$-subgroups.
\end{theo}
A significant portion of the paper is taken up with showing that if
$G$ is an almost simple group satisfying $P(p,q)$, then
$\gcd(|G|,q)=1$, this uses CFSG.

Interestingly an example where a given prime almost does not occur is given by:
\begin{theo}\cite{MR2502216} Let $G$ be a finite group having exactly one
conjugacy class of size a multiple of a prime $p$. Then one of the following
holds: \\(i) $G$ is a Frobenius group with Frobenius complement of order $2$
and Frobenius kernel of order divisible by $p$;\\ (ii) $G$ is a doubly
transitive Frobenius group whose Frobenius complement has a nontrivial
central Sylow $p$-subgroup;\\ (iii) $p$ is odd, $G=KH$ where $K= \Fit(G)$
 is a $q$-group, $q$ prime, $H= C_G(P)$ for a Sylow
$p$-subgroup $P$ of $G$, ${K\cap H}=Z(K)$ and $G/Z(K)$ is a doubly transitive
Frobenius group.\end{theo}
In the paper there is a slightly more detailed version which gives
if and only if conditions.

\subsection{Solubility}\label{sol}

Recognising solubility is clearly an interesting problem. In
1990 Chillag \& Herzog proved the following:
\begin{theo}\label{four} \cite{MR1055001} Suppose $4$ does not divide $\ind_G(x)$
for all $x \in G$. Then $G$ is soluble.
\end{theo}
A proof avoiding CFSG
is given in \cite{MR1633180} and \cite{MR1705872}.
Chillag \& Herzog also considered those groups all of whose indices are
square-free; these groups were also studied by Cossey \& Wang \cite{MR1705872}.
\begin{theo} \cite{MR1705872}
Suppose that all indices of the group $G$ are
square-free. Then $G$ is supersoluble and both $|G/\Fit(G)|$ and
$G'$ are cyclic groups with square-free orders.
The class of $\Fit(G)$ is at most 2 and $G$ is metabelian.
\end{theo}

Using the results of \cite{MR636194} quoted
in Section~\ref{CFSG}  Li strengthens these results as follows:
\begin{theo}\cite{MR1741423} Let $G$ be a group and let $p$ be the smallest
prime dividing the order of $G$. Assume that $p^2$ does not divide the
index of any element of $q$-power order, for $q$ any prime not equal to $p$.
Then $G$ is $p$-nilpotent. In particular, $G$ is soluble.
\end{theo}

\begin{theo} \cite{MR1741423}
Suppose that all indices of elements of prime-power order in the group $G$ are
square-free. Then $G$ is supersoluble, the
derived length of $G$ is bounded by 3, $G/\Fit(G)$ is a direct
product of elementary abelian groups and $|\Fit(G)'|$ is a
square-free number.
\end{theo}

In \cite{MR2513573} the authors consider what happens if only
$p'$-elements have indices not divisible by $p^2$. They show that in
this case the highest power of $p$ which can divide the order of any
chief factor is $p$.

Note, if $p$ is the smallest prime that divides the order of a group
then for any other prime $q$ dividing the order of the group $q$
does not divide $p-1$. Using this observation Cossey \& Wang
~\cite{MR1705872}  considered groups $G$ for which there is a prime
$p$ dividing the order of $G$ so that if $q$ divides the order of
$G$ then $q$ does not divide $p-1$.  They prove a result giving the
structure of such groups when the index of no element of $G$ is
divisible by $p^2$. A variation on this  is given by Liu, Wang and
Wei ~\cite{MR2111851}. We note that with these conditions if $p=2$
then $G$ is soluble by Theorem~\ref{four} and if $p>2$ then $G$ has
odd order.

If the prime $p$ divides the order of the group and there exists primes which do not divide $p-1$ we can say something slightly more general:
\begin{lem}
Let $G$ be a soluble group and let $p$ be a prime. Assume that all elements
of $p'$-order have indices not divisible by $p^2$.
Let $\Pi=\{q:q \mbox{ a prime}, q \neq p \mbox{ and } q
\mbox{ not dividing } p-1\}$.
Then $G/\OG_{p'}(G)$ is a $\Pi'$-group.
\end{lem}
{\bf Proof.} Let $G$ be a minimal counter-example.
If $\Pi$ is empty there is nothing to prove.
Note that this hypothesis goes over to normal subgroups and quotients.
Clearly if $\OG_{p'}(G)\neq 1$ we obtain the result by factoring
out $\OG_{p'}(G)$. So we can assume  $\OG_{p'}(G)=1$.

Assume that $G$ has a proper normal subgroup $N$, from the above comment
$\OG_{p'}(N)=1$. So $N$ is a $\Pi'$-group. If $N$ is a maximal
normal subgroup it has prime index, say $q$, then $q\in \Pi$. If there were two distinct
maximal normal subgroups then $G$ would be a $\Pi'$-group and so there is nothing to prove.
So $N$ is the unique maximal normal subgroup of $G$.

 Let $M$ be a minimal normal subgroup. Let $x$ be an element of order $q$, which is not in $N$ as $q\in \Pi$ and $N$ is a $\Pi'$-group, and consider the action of
 $x$ on $M$. Since $M$ is normal $[M:\Cent_M(x)]=1\mbox{ or } p$. Since $q \in \Pi$ the action of $x$ is trivial
 on $M/\Cent_M(x)$, so $x$ centralises $M$. Since $x \notin N$ and $x\in \Cent_G(M)$ it follows that $M\leq Z(G)$. But $\OG_{p'}(G)=1$
 so $\OG_{p'}(G/M)=1$ and $G/M$ is a $\Pi'$-group which contradicts the
assertion that $[G:N]=q$.$\Box$\\[1ex]
Note that if  $G$ is a  $\{p\}\cup \Pi$-group then $G$ has a normal
$p$-complement.

Recall, an $A$-group is a group with abelian Sylow subgroups. Camina
\& Camina proved the following result:
\begin{theo}\cite{MR2228630} Let $G$ be an $A$-group and suppose
$2^a$ is the highest power of 2 to divide
an index of $G$. Then, if there exists an element $x \in G$ with
$\ind_G(x) = 2^a$, it follows that $G$ is soluble.
\end{theo}
An easy adaptation of the proof, shows that the result holds if
 instead of requiring our group to be an $A$-group, we
just require the Sylow 2-subgroup to be abelian.
However, if the Sylow 2-subgroup
is not abelian the result is false. For example, the simple group
${\rm PSL}(2,7)$ has a permutation representation of degree 8.
Let $V$ be a permutation module of degree 8 over a finite field of characteristic
5. By considering the extension of $V$ by
${\rm PSL}(2,7)$ it can be seen that $V$ has elements of
index 8.\\[2ex]
Marshall considered soluble $A$-groups. Let the conjugate rank of
$G$, $\mbox{crk}(G)$, be greater than one. She proved that there
exists a function $g: \mathbb{Z}^+ \mapsto \mathbb{Z}^+$ such that
the derived length of a soluble $A$-group $G$ is bounded by
$g(\mbox{crk}(G)+1)$ \cite{MR1411077}. In 1997 Keller gave an explicit function for the bound
which is stronger than logarithmic,
 see his survey article \cite{MR2051537}. This can be seen as a contribution
 towards the following question.

\begin{ques} Is it possible to bound the derived length of a
soluble group by its conjugate type rank?\end{ques}

We note the following interesting result by Keller:- \begin{theo}
\cite{MR2187399} Let $G$ be a finite soluble group then the derived length
of $G / \Fit(G)$ is bounded by $ 24\log_2(\mbox{crk}(G)+1)+364$.
\end{theo}
We finish this section with a question related to the comments in
Section~\ref{CHA}.
\begin{ques} (i) If we know all conjugacy class sizes including
multiplicities can we recognise solubility?\\
(ii) If we know all conjugacy class sizes can we recognise
solubility?\end{ques}

Clearly (ii) is stronger than (i) but we have no idea what the
answer might be. One can also look at recognising other classes like
supersolubility.

\section{Conjugate Type Vectors}\label{CTP}

\subsection{It{\^o}}

Recall, the conjugate type vector is a list of the distinct
conjugacy class sizes in descending order,
 and the conjugate rank is the number of entries not equal to $1$.
It{\^o} proved the following:
\begin{theo}  Let $G$ be a finite group,\\
(i) \cite{MR0061597} with conjugate type vector $(n,1)$. Then $n = p^a$ for some
prime $p$ and $G$ is nilpotent. More exactly,
$G$ is a direct product of a $p$-group and
an abelian $p'$-group.  The $p$-group, $P$, has an abelian normal
subgroup $A$, such that $P/A$ has exponent $p$.\\
(ii) \cite{MR0263910} with conjugate type vector $(n_1, n_2, 1)$. Then $G$ is soluble.
\end{theo}

Recently Ishikawa \cite{MR1910937} proved that a $p$-group of conjugate rank
$1$ has nilpotency class at most $3$. Almost immediately a number
of authors generalised the result, some to Lie Algebras,
\cite{MR1953720, MR2085724, MR2267792, MR2390495,MR2470541,Mann}. The
generalisations involve the subgroup $M(G)$ of finite group
$G$, where $M(G)$ is defined to be the subgroup generated by elements
whose indices are $1$ and $m$ and $m$ is the smallest
non-trivial index. Then Isaacs  proved the following:
\begin{theo}\cite{MR2390495} Let $G$ be a finite group which contains a normal
abelian subgroup $A$ with $\Cent_G(A)=A$. Then $M(G)$ is nilpotent of
class at most $3$.
\end{theo}

In 2009 Guo, Zhao \& Shum proved the following generalisation of the rank one case:
\begin{theo}\cite{MR2530762} Let $N$ be a
$p$-soluble normal subgroup of a group $G$ such that $N$ contains a
noncentral Sylow $r$-subgroup, ($r \neq p$), $R$ of $G$. If $|x^G|=1$ or $m$
for every $p'$-element $x$ of $N$ whose order is divisible by at most
two distinct primes, then the $p$-complements of $N$ are nilpotent.
\end{theo}
It{\^o}'s $1952$ result follows as a corollary. Zhao and Guo give another generalisation as follows:
\begin{theo}\cite{MR2529448} Let $G$ be a finite group with a non-central Sylow
$r$-subgroup $R$ and $N$ a normal subgroup of $G$ containing $R$. If $|x^G| =1$
or $m$ for every element $x$ of $N$, then $N$ is nilpotent.
\end{theo}

In 1996 Li proved the following extension of It{\^o}'s rank 1 result,
again using the
argument quoted in Section~\ref{CFSG}:
\begin{theo} \cite{MR1399825} Let $G$ be a finite group and let $m$ be a natural number.
Assume that if $x \in G$ has prime-power order then $x$ has index $1$ or $m$. Then $G$ is
soluble.\end{theo}
Another variation is given by considering the $p'$-conjugate type vector.
\begin{defn} The $p'$-conjugate type vector of a group $G$ is the list
of distinct indices of $p'$-elements in descending order.
\end{defn}

Suppose $G$ is a group with $p'$-conjugate type vector
$(m,1)$. Three authors, \cite{MR1962969, MR2505347}, have shown that
 $m =p^aq^b$ for primes $p \neq q$ and if $a$ and $b$ are both strictly
greater than 0 then $G =PQ\times A$ where $P$ is a Sylow
$p$-subgroup of $G$, $Q$ is a Sylow $q$-subgroup of
$G$ and $A$ is in the centre of $G$. We comment that the result follows from the results
in \cite{MR0316563} and \cite{MR0346054}. The first two authors also considered the case where the $p'$-conjugate type vector
is $(m,n,1)$ where $\gcd(m,n)=1$, \cite{MR2099353}.\\[1ex]
\subsection{Conjugate rank $2$}

In 1974, A.R.~Camina gave a new proof of the solubility of groups of conjugate rank 2, along with more details on the structure of such groups.
 This depended on work
of Schmidt~\cite{MR0313388} and Rebmann~\cite{MR0291275}
which looked at the lattice of centralizers.
\begin{defn} A group $G$ is an $F$-group if given any
pair $x,\,y$ with $x,\,y \notin Z(G)$ then $\Cent_G(x)\not<
\Cent_G (y)$. \end{defn}

In these papers the authors completely classify $F$-groups. It is
not possible to express the condition within the framework of
indices. However the condition that for no pair of indices is one
divisible by the other implies that the group is an $F$-group.

A.R. Camina's theorem says what happens if this does not occur:
\begin{theo} \cite{MR0346054} If $G$ has conjugate rank 2 and is not an $F$-group then
$G$ is a direct product of an abelian group and a group whose order
is divisible by only two primes (or $\rho^*(G)=2$).\end{theo}

This has recently been improved by Dolfi \& Jabara. Their proof is
independent of \cite{MR0346054} but uses \cite{MR0263910}:
\begin{theo}\cite{DoJa} A finite group $G$ has conjugate rank $2$ if and
only if, up to an abelian direct factor, either\\
(1) $G$ is a $p$-group for some prime $p$ or\\
(2) $G=KL$, with $K\unlhd G,\,\gcd(|K|,|L|)=1$ and one of the following
occurs\\
\indent(a) both $K$ and $L$ are abelian, $Z(G)<L$ and $G$ is a
quasi-Frobenius group,\\
\indent(b) $K$ is abelian, $L$ is a non-abelian $p$-group,
for some prime $p$ and
${\rm O}_p(G)$ is an abelian subgroup of index $p$ in $L$ and
$G/\mathrm{O}_p(G)$ is a Frobenius group or\\
\indent(c)  $K$ is a $p$-group of conjugate rank $1$ for some prime $p$,
$L$ is abelian, $Z(K)=Z(G)\cap K$ and $G$ is quasi-Frobenius.
\end{theo}

We note that the results in \cite{MR2542207} can be deduced from the
results of  Dolfi \& Jabara.
This completes the classification of groups of conjugate rank $2$.

\subsection{Conjugate rank larger than 2}
It{\^o} went on in the early 1970's to consider groups of low
conjugate rank, $3$, $4$
and $5$ with special reference to the simple groups which can occur,
\cite{MR0289625, MR0289636, MR0316549, MR0354845}.
Recently there has been some effort to look at the case of
conjugate rank $3$. The first papers were those connected with proving
nilpotence, see Subsection~\ref{nilpotency}. Beltr{\'a}n and Felipe
\cite{MR2379085} looked at the
structure of soluble groups with conjugate type vector $(mk,m,n,1)$ where
$\gcd(m,n)=1$ and $k$ divides $n$.
More recently Camina \& Camina
have shown that if the conjugate rank is larger than $2$ and there are two
coprime indices amongst any three, then $G$ is soluble. Amongst the results are the following:
\begin{cor}~\cite{MR2532697} Let
$G$ be a finite group with trivial centre and with at most three
distinct conjugacy class sizes greater than 1. Then $G$ is either
soluble or PSL$(2,2^a)$.\end{cor}
\begin{cor}~\cite{MR2532697} Let $G$ be a finite $A$-group with at most
three distinct conjugacy class sizes greater than 1. Then $G$ is
either soluble or PSL$(2,2^a)$.\end{cor}

So we ask \begin{ques}
 Can the rank $3$ groups be classified, especially the non-soluble ones?
\end{ques}

We note the following theorem due to Bianchi, Gillio and Casolo, which followed
earlier work \cite{MR1143160}, \cite{MR1682749}:
\begin{theo}\cite{MR1876222} Suppose $G$ is a group with conjugate type
vector $(m,n \ldots)$ where $m$ and $n$ are coprime. Suppose $x,y \in G$
with ${\rm Ind}_G(x) = n$ and ${\rm Ind}_G(y) =m$ and
let $N ={\rm C}_G(x)$ and $H={\rm C}_G(y)$. Then $N$ and $H$ are abelian
and $G$ is quasi-Frobenius with kernel $N/Z(G)$ and complement $H/Z(G)$.
\end{theo}
Thus the conjugate type vector is of the form $(m,n,1)$ (compare with
Theorem~\ref{isolated}).

\subsection{Nilpotency}\label{nilpotency}

In 1972 A.R. Camina proved that a finite group with conjugate type
vector $(q^bp^a, q^b, p^a, 1)$, where $p$ and $q$ are primes, is
nilpotent \cite{MR0294481}. Recently Beltr{\'a}n and Felipe have
proved that if $G$ is has conjugate type vector $(nm, m, n, 1)$ with
$n$ and $m$ coprime integers, then $G$ is nilpotent and $n$ and $m$
are prime powers \cite{MR2192606,MR2464108,MR2272717}.

Beltr{\'a}n and Felipe~\cite{MR2295080} have also shown for
$G$ a finite $p$-soluble group and $m$ a positive integer not
divisible by $p$ that if the set of conjugacy class sizes of all
$p^\prime$-elements of $G$ is $\{1,m, p^a,mp^a\}$, then $G$ is nilpotent
and $m$ is a prime power.

Camina's result led to the question whether you could identify a nilpotent group
from its conjugate type vector. More precisely if $G$ and $H$ have the same
conjugate type vector and $H$ is nilpotent, does it follow that $G$ is
nilpotent? Notice that a nilpotent group satisfies the following property:
if $m$ and $n$ are coprime integers, there are elements of index $m$ and of
index $n$ if and only if there is an element of index $mn$. This follows
from the fact that a nilpotent group is a direct product of its Sylow
$p$-subgroups. We note the following nice result due to
Cossey and Hawkes:\label{CHA}
\begin{theo} \cite{MR1641677} Let $p$ be a prime and $\mathcal{S}$ a finite set of
$p$-powers containing 1. Then there exists a $p$-group $P$ of class
2 with the property that $\mathcal{S}$ is the conjugate type vector
of $P$ (ordered appropriately).
\end{theo}
(Note it is certainly not true that arbitrary sets of numbers can be
conjugate type vectors, as the results on graphs of Section 5
indicate.) Using Cossey and Hawkes result it follows that
recognising nilpotency is equivalent to asking whether a group whose
indices satisfy the property indicated above is nilpotent.

Although in certain cases you can recognise nilpotency, for example
if all the conjugacy classes are square-free \cite{MR1633180}, or if
the group is a metabelian $A$-group \cite{MR2228630}, in general it
is not true. The smallest example of a group who shares its
conjugate type vector with a nilpotent group, but is not itself
nilpotent, has order 160 and conjugate type vector $(20, 10, 5, 4,
2, 1)$. An infinite family of such examples is given in
\cite{MR2228630}. A number of questions are posed in the paper.

\begin{ques} Suppose $G$ and $H$ are finite groups with $H$ nilpotent,
further suppose $G$ and $H$ have the same conjugate type vector.\\
(i) Is it true that $G$ must be soluble? \\
(ii) If $G$ is not nilpotent, does $G$ have a centre?\\
(iii) Suppose $G$ is an $A$-group, then must $G$ be nilpotent?
\end{ques}

Note, if you have the additional information of the number of
conjugacy classes of each size, then you can recognise nilpotency as
Cossey, Hawkes \& Mann have proved:
\begin{theo} \cite{MR1157262} Let $G$ and $H$ be finite groups with
$H$ nilpotent. Let $S_H$ be the multiset of conjugacy class sizes of
elements in $H$ and define $S_G$ similarly. Suppose $S_H = S_G$,
then $G$ is nilpotent.
\end{theo}
These ideas have been extended by Mattarei \cite{MR2213690}.

\section{Graphs}\label{graphs}
In an earlier version of this survey we included a lengthy chapter
on graphs associated to conjugacy class sizes. However, since then
 Lewis has published an excellent survey concerned with graphs
associated to character degrees and conjugacy class sizes
\cite{MR2397031}. Thus we will just briefly introduce the graphs and
mention a few of the more recent results, while still aiming to
include a large bibliography.

Let $X$ be a set of positive integers. We associate two graphs
to $X$, the {\it prime vertex graph}  and the {\it
common divisor graph}.
\begin{defn} (i) The common divisor graph of $X$ has
vertex set $X^* = X \setminus 1$
($X$ may or may not contain the element 1) and an edge
between $a, b \in X^*$ if $a$ and $b$ are not coprime. We denote
the common divisor graph of $X$ by $\Gamma(X)$.\\
(ii) The prime vertex graph has vertex
set $\rho(X) = \bigcup_{x \in X} \pi(x)$ where $\pi(x)$ denotes the prime
divisors of $x$. There is an edge between $p, q \in \rho(X)$ if $pq$ divides
$x$ for some $x \in X$. The prime vertex graph is denoted by $\Delta(X)$.
\end{defn}
The connection between these two graphs has
been clarified by the recent paper \cite{IrPr} which defines a bipartite
graph $B(X)$.
\begin{defn} The vertex set of $B(X)$ is given by
the disjoint union of the vertex set of $\Gamma(X)$ and the vertex set
of $\Delta(X)$, i.e. $X^* \cup \rho(X)$. There is an edge
between $p \in \rho(X)$ and $x \in X^*$
if $p$ divides $x$, i.e. if $p \in \pi(x)$.
\end{defn}
Consideration of $B(X)$ makes
it clear that the number of connected components of $\Gamma(X)$ is equal
to the number of connected components of $\Delta(X)$. Furthermore the
diameter of a connected component of $\Gamma(X)$ differs by at most one
from the diameter of the equivalent connected component of $\Delta(X)$.

There are two common choices for the set $X$:
\begin{defn}(i) ${\rm cs}(G)$
the set of conjugacy class sizes of $G$ (equivalently the set of
indices of elements of $G$), and (ii) ${\rm cd}(G)$ the set of
degrees of irreducible characters of $G$. \end{defn}

Much has been written on
both these cases but we will concentrate on the case when $X =
{\rm cs}(G)$. In 1981 Kazarin wrote a paper on isolated
conjugacy classes \cite{MR636915}. A group $G$ has isolated
conjugacy classes if there exist elements $x,y \in G$ with coprime
indices such that every element of $G$ has index coprime to either
${\ind}_G(x)$ or ${\ind}_G(y)$\label{isol}. Kazarin classified all
groups with isolated conjugacy classes. He therefore classified
all groups such that either $\Gamma({\rm cs}(G))$ has more than
one component or $\Gamma({\rm cs}(G))$ is connected and has
diameter at least $3$:
\begin{theo}\label{isolated}\cite{MR636915}
Let $G$ be a group with isolated conjugacy classes. Let $x$ and $y$ be representatives of the two
isolated conjugacy classes with $\ind_G(x)=m_1$ and $\ind_G(y)=n_1$.
Then $|G|=mnr$ where $r$ is coprime to both $m$ and $n$, the only primes
which divide $m$, respectively $n$, are those which divide $m_1$,
respectively $n_1$. Further $G=R \times H$ where $|R|=r$ and $H$ is a
quasi-Frobenius group where the pre-image of both the kernel and complement are abelian.
\end{theo}

It follows that
$\Gamma({\rm cs}(G))$ has at most two connected components. If this is
the case we can take (without loss of generality) the conjugate type
vector to be $(m_1, n_1, 1)$. Furthermore, if
$\Gamma({\rm cs}(G))$ is connected its diameter is at most 3. This
has also been proved using different techniques in
\cite{MR1099007} and \cite{MR1209249}.

In \cite{MR1099007} the authors note that if $G$ is a nonabelian
simple group then $\Gamma({\rm cs}(G))$ is complete, this follows
from work of Fisman \& Arad \cite{MR892909}.
Another nice result is due to
Puglisi and Spezia. They prove that if
$\Gamma({\rm cs}(H))$ is a complete graph for every subgroup $H$ of
$G$, then $G$ is soluble \cite{MR1683129}.

It{\^o} proved an early result about $\Delta({\rm cs}(G))$ : suppose $p$ and
$q$ are two distinct non-adjacent vertices in $\Delta({\rm cs}(G))$
then $G$ is either $p$-nilpotent or $q$-nilpotent \cite[Proposition 5.1]{MR0061597}. Dolfi extended this
result for the case when $G$ is soluble concluding that in this case both the
$p$-Sylow and $q$-Sylow subgroups are abelian \cite{MR1337169}.
In this paper Dolfi also proves
the analogous structural results for $\Delta({\rm cs}(G))$: namely,
if $\Delta({\rm cs}(G))$
is not connected then it has exactly two connected components and they are
complete graphs, and if $\Delta({\rm cs}(G))$ is connected it has diameter at
most 3. These results are also contained in \cite{MR1276135}.
Alfandary determines further structural results for $\Delta({\rm cs}(G))$
when $G$ is soluble in a follow up paper \cite{MR1351622}.

Casolo and Dolfi~\cite{CaDo} have characterised the groups for which
$\Delta({\rm cs}(G))$ has diameter 3 \cite{MR1367161}.
In the same paper they show that
if a group is not soluble then $\Delta({\rm cs}(G))$ is connected and
has diameter at most 2.
Another nice result of Dolfi is to show that
given 3 distinct vertices in $\Delta({\rm cs}(G))$ then at least two are
connected by an edge \cite{MR2253660}.

A recent paper gives further results on the bipartite graph $B({\rm
cs}(G))$ \cite{MR2518171}.

Slightly confusingly, another graph has been associated to a finite
group $G$. In this case the vertices are given by the set of non-central
conjugacy classes and two vertices $C$ and $D$ are joined if $|C|$ and
$|D|$ share a common divisor \cite{MR1099007}.
We shall call this graph $\hat{\Gamma}(G)$
and note that it shares many properties with $\Gamma({\rm cs}(G))$.

We note that the following conjecture can be viewed as considering groups $G$
for which $\Gamma({\rm cs}(G)) = \hat{\Gamma}(G)$.\\[1ex]
{\bf $S_3$ conjecture}: Any finite group in which distinct conjugacy classes have distinct sizes is isomorphic to $S_3$.\\[1ex]
This conjecture has been verified for soluble groups by Zhang
\cite{MR1302859}, and Kn{\"o}rr, Lempken and Thielcke
\cite{MR1348305}, independently. More recently Arad, Muzychuk and
Oliver have studied the case of insoluble groups \cite{MR2089251}.

One type of question that is posed when considering these graphs is
how the structure of the graph determines the structure of the
group? For example the problem of classifying all groups
$G$ such that $\hat{\Gamma}(G)$
has no subgraph $K_n$ (where $K_n$ is the complete graph with $n$ vertices)
has been considered in \cite{MR2161799} for  $n = 4,5$ and in \cite{MR1981424}
for $n=3$. In \cite{MR2532697}
groups $G$ such that $\Gamma({\rm cs} (G))$ has no triangles have been
considered, such groups have conjugate rank at most 3 and are soluble.
This yields the following question.
\begin{ques} Let $G$ be a finite group and $n$ a natural number. Is it
true that if $\Gamma({\rm cs}(G))$ has no subgraphs isomorphic to $K_n$
then there is a function of $n$ which bounds the conjugate rank
of $G$?
\end{ques}
Note that this question can also be asked without the language of graphs.
The condition translates to requiring that given any set of
$n$ distinct indices then there are two which are coprime.

Variants of these graphs have been introduced. For example
Beltr{\'a}n \& Felipe have considered a version of $\hat{\Gamma}(G)$
where the vertices are restricted to the $p'$-conjugacy classes,
that is the set of $x^G$ where $x$ is a $p'$-element. They consider
the case when $G$ is $p$-soluble \cite{MR1941928, MR1995542,
MR2097475} and have summarised their results in a nicely written
survey article \cite{MR2328161}.

Alternatively, Qian and Wang \cite{MR2510974} have considered the
conjugacy class sizes
of $p$-singular elements, that is elements whose order is divisible
by $p$. They denote the set of $p$-singular elements
by $G_p$ and consider the graph $\Gamma({\rm cs}(G_p))$. Noting that if
$p$ divides $|Z(G)|$ then $\Gamma({\rm cs}(G_p)) = \Gamma({\rm cs}(G))$
they prove
that if $p$ divides $|G|$ but not $|Z(G)|$ then $\Gamma({\rm cs}(G_p))$ is
connected with diameter at most 3. In the paper the authors also
consider groups for which $p$ does not divide $m$ for any $m \in {\rm cs}(G_p)$,
this leads to a classification of all finite groups for which every
conjugacy class size coincides with a Hall number.

Beltr{\'a}n has introduced
the $A$-invariant conjugacy graph \cite{MR1968426}.
Let $A$ and $G$ be finite groups and
suppose that $A$ acts on $G$ by automorphisms. Then $A$ acts on the set of
conjugacy classes of $G$. The $A$-invariant conjugacy graph $\Gamma_A(G)$
has vertices the non-central $A$-invariant conjugacy classes of $G$ and two
vertices are connected by an edge if their cardinalities are not coprime.
Consideration of the case when $A$ acts trivially gives that $\Gamma_A(G)$ is
a generalisation of $\hat{\Gamma}(G)$. Beltr{\'a}n notes that the proof
that $\hat{\Gamma}(G)$ has at most two connected components given in
\cite{MR1099007}
translates to this more general setting. He then considers the disconnected
case and proves the following theorem:
\begin{theo} \cite{MR1968426} Suppose that a group $A$ acts coprimely on a group
$G$ and that $\Gamma_A(G)$ has exactly two connected components. Then $G$ is
solvable.
\end{theo}
It is not known whether the result holds when $A$ and $G$ do not have coprime
orders.

In a different direction Isaacs and Praeger introduced a
generalisation of $\Gamma({\rm cs}(G))$ known as the IP-graph
\cite{MR1231209}.
\begin{defn} Let $G$ be a group acting transitively on a set
$\Omega$, and let $D$ denote the set of subdegrees of $(G,
\Omega)$, that is, the cardinalities of the orbits of the action
of a point stabilizer $G_{\alpha}$ on $\Omega$. Suppose the
subdegrees are finite, then the IP-graph of $(G, \Omega)$ is the
common divisor graph of $D$.
\end{defn}
That this is a generalisation of
$\Gamma({\rm cs}(G))$ can be seen as follows. Let $G$ be a group
and ${\rm Inn}(G)$ the inner automorphisms of $G$. Let $H$ be the
semidirect product $G \rtimes{\rm Inn}(G)$, then $H$ acts
transitively on $G$ by sending $x \in G$ to $(xg)^{\sigma}$ where
$g \sigma \in G \rtimes {\rm Inn}(G)$. Clearly the orbits of
$H_1$, the stabilizer of the identity, are the conjugacy classes.
The authors prove that the IP-graph of $(G, \Omega)$ has at most
two connected components and that the diameter of a connected
component is bounded by 4. However they know of no example of
diameter 4, suggesting perhaps that 3 is the correct upper bound.
\begin{ques} Does there exist a group $G$ and a set $\Omega$
such that the diameter of the IP-graph of $(G, \Omega)$ is 4?
\end{ques}
Kaplan has studied the case when the IP-graph is disconnected
\cite{MR1458805}, \cite{MR1710739}.  Neumann introduced the VIP
graph, a variant of the IP graph which does not restrict itself
to the case where all subdegrees are finite \cite{MR1207949}. More
recently the IP graph has been generalised to the setting of naturally
valenced schemes~\cite{MR2430305}. This work was extended by
 Xu~\cite{MR2504489}.

We would like to introduce a graph using the notion of divisibility.
Let $X$ be a set of positive integers, then $D(X)$ the {\it divisibility
graph} is a directed graph. The vertex set of $D(X)$ is given by
$X^{\ast}$ and there is an edge connecting $(a,b)$ with
$a, b \in X^{\ast}$ whenever $a$ divides $b$. We are interested in the
properties of this graph when $X = {\rm cs}(G)$ for a finite group $G$.
Note that if $D({\rm cs}(G))$ has no edges then $G$ is an $F$-group.
\begin{ques} How many components can $D({\rm cs}(G))$ have?
\end{ques}
It is worth pointing out that whilst many results can be interpreted
in the language of graph theory there are many interesting problems
that have no such simple description.
\section{The number of conjugacy classes}\label{noclasses}
Given a group $G$ of order $n$ with $k$ conjugacy classes what can
be said about the relation between $n$ and $k$? It is trivial to
see that $k\leq n$ but can anything be said in the opposite
direction? The first to bound $n$ in terms of $k$
was Landau in 1903 \cite{MR1511192}. He used a number theoretic
approach to the class equation. Brauer \cite{MR0178056} was the
first to give an explicit bound using Landau's method. He asked, in
Problem 3, whether better methods could be found. A similar
approach was taken by Newman in \cite{MR0225870} where he improved
Landau's result. This gave very general bounds of
exponential form. He proved:
\begin{theo}\cite{MR0225870} Let $G$ be a finite group of order $n$ with $k$
conjugacy classes. Then
$$k\geq \frac{\log\log n}{\log 4}.$$
\end{theo}

Also there have been papers which give complete descriptions for
small $k$, the earliest example being Miller in 1919
\cite{MR1501126}. The objective here is to try to classify the
isomorphism classes with a given number of conjugacy classes. A
number of authors over many years have classified groups with
few conjugacy classes
\cite{MR0011299, MR0223447, MR0233883, MR804489, MR880291}.
In particular, Vera L{\'o}pez \& Vera L{\'o}pez
\cite{MR804489,MR880291} examine groups with $13$ and
$14$ conjugacy classes and give lists of such groups. As far as
the authors know this is the largest value of $k$ for which this
has been attempted.

Cartwright in \cite{MR906145} considered soluble groups. He
proves:
\begin{theo} There exist positive constants $a$ and $b$ so that a
soluble group of order $n$ has at least $a(\log n)^b$ conjugacy
classes.
\end{theo}
There is an estimate for $b$ of $0.00347$. But Pyber
\cite{MR1182481} has proved the following theorem:
\begin{theo} Let $G$ be a finite group with $k$ conjugacy classes.
Then
$$k\geq \epsilon\log n/(\log \log
n)^8,$$ for some fixed $\epsilon$.  \end{theo}
 This could be considered an answer to Brauer's question. We note that Keller has reduced the power in the denominator from 8 to 7
\cite{keller}. Also, for nilpotent groups Jaikin-Zapirain has proved a stronger result \cite{JZ}.

In 1997 Liebeck \& Pyber \cite{MR1489911} wrote a paper in which
they gave upper bounds for the number of conjugacy classes  for
various classes of groups. Their main result is the following:
\begin{theo}\cite[Theorem 1]{MR1489911} Let $G$ be a finite simple
group of Lie type over a finite field of order $q$. Let $G$ have
rank $\ell$ and assume that $G$ has $k$ conjugacy classes. Then
$$k\leq (6q)^{\ell}.$$\end{theo}
Note that this is the``untwisted" rank.

As a consequence of this they prove
\begin{theo}\cite[Theorem 2]{MR1489911} Let $G$ be a subgroup of the symmetric group
of degree $n$. Then if $G$ has $k$ conjugacy classes
$$k\leq 2^{n-1}.$$\end{theo}
 This solved a conjecture of Kov\'acs \&
 Robinson~\cite{MR1244923}. They had found a bound of $5^{n-1}$,
 without the classification of finite simple groups. Mar{\'o}ti
 has reduced the bound to $k\leq 3^{(n-1)/2}$ \cite{MR2137971}.

An interesting variation has been proved by Jaikin-Zapirain:
\begin{theo} \cite{MR2137970} There exists a function $f(r)$ such that if
 $G$ is a soluble group with at most $r$ conjugacy classes of size $k$ for any $k$ then $|G|\leq f(r)$.
\end{theo}

\section{Comparisons with Character Theory}\label{char}
\subsection{Introduction}
The most obvious connection with character theory is that the number
of conjugacy classes
is the same as the number of irreducible characters. A number of authors,
including those of
this article, have seen a connection between character degrees and conjugacy class sizes
and have searched for analogous results.
As we commented in the Introduction, if we know
the multiplication constants then we have a complete connection between the
two sets of data.
However if we only know the
degrees and the indices the link seems less clear.

\subsection{The differences}

Perhaps the first obvious difference is that the order of the group
is given by the sum of the sizes of the conjugacy classes, but the
sum of the squares of the degrees of the irreducible characters.
This might help to explain the following dichotomies. As previously
mentioned in Section 2, if $p$ is coprime to all indices of elements
of $G$ then the Sylow $p$-subgroup of $G$ is an abelian direct
factor of $G$. However, It{\^o} proved that if $G$ has a normal
abelian Sylow $p$-subgroup then $p$ is coprime to all character
degrees of $G$. Also, in Section \ref{sol} we noted that if all
indices are square-free then $G$ is soluble. However, if all
irreducible character degrees are square-free $G$ need not be
soluble, the smallest example we know is Alt$(7)$, see for example
\cite{MR827219}. Another example is given by the two different
conclusions drawn when $\{1, p^a, q^b, p^aq^b\}$ is either the set
of character degrees or the set of conjugacy class sizes and $p$ and
$q$ are distinct primes. In the conjugacy class case we can conclude
that $G = P \times Q$ where $P$ is the Sylow $p$-subgroup and $Q$
the Sylow $q$-subgroup \cite{MR0294481}. This conclusion does not
hold in the character case \cite{MR1637284}.

One might also ask whether there are connections between the sizes
of ${\rm cs}(G)$ and ${\rm cd}(G)$? The answer is no. In an elegant construction Fern{\'a}ndez-Alcober \& Moret{\'o} prove the following theorem:-
\begin{theo}\cite{MR1844993} Given any two integers $r$ and $s$ greater than 1 there exists a $p$-group $G$ of nilpotency class $2$ such that $|cd(G)|=r$ and $|cs(G)|=s$. \end{theo}

\subsection{The similarities}

In Section 5 we defined the common divisor graph $\Gamma(X)$ and the
prime vertex graph $\Delta(X)$ for a set of natural numbers $X$. We
focussed on the case when $X$ is the set of indices of a finite
group $G$, denoted by ${\rm cs}(G)$, but commented that much work
has been done for the case when $X$ is the set of degrees of
irreducible characters, ${\rm cd}(G)$, in fact this case was
considered first (see \cite{MR2397031} for a recent survey of
results). What is remarkable is the similarity of the graphs for the
two different choices of $X$. In particular, if $G$ is soluble then
in all cases (both choices of $X$ and graph) the graph has at most
two connected components and if the graph is connected the diameter
of the graph is at most three. However, if $G$ is not soluble then
it is possible for $\Delta({\rm cd}(G))$, and hence also
$\Gamma({\rm cd}(G))$, to have three connected components.

Recently Casolo and Dolfi have proved that $\Delta({\rm cd}(G))$
is a subgraph of $\Delta({\rm cs}(G))$.
It is easy to see that each prime number that divides an
irreducible character degree of $G$ must also divide some class size
of $G$. Casolo and Dolfi~\cite{CaDo} prove that for
distinct primes $p$ and $q$ if $pq$ divides the degree of some
irreducible character of $G$, then it also divides the size of some
conjugacy class of $G$. Dolfi had
previously proved this result for soluble groups \cite{MR1337168}.

Authors have considered the implications of the arithmetical data of
$\mbox{cd}(G)$ on the structure of the group. For example
 in \cite{MR1637284} Lewis proves that for $p, q$ and $r$
distinct primes, if $\mbox{cd}(G) = \{1,p,q,r, pq, pr\}$ then $G = A
\times B$ where $\mbox{cd}(A) = \{1,p\}$ and
$\mbox{cd}(B)=\{1,q,r\}$. For $s$ another distinct prime he also
proves that if $\mbox{cd}(G)=\{1,p,q,r,s,pr,ps,qr,qs\}$ then $G =
A\times B$ with $\mbox{cd}(A)= \{1,p,q\}$ and
$\mbox{cd}(B)=\{1,r,s\}$. The equivalent results where character
degrees are replaced by conjugacy class sizes are proved in
\cite{MR1800744}.

In 2006 Isaacs, Keller, Meierfrankenfeld \& Moret{\'o} showed that
if $G$ is soluble, $p$ a prime divisor of $|G|$ and $\chi$ a primitive
irreducible character of $G$. Then the $p$-part
of $\chi(1)$ divides the $p$-part of $(\mbox{Ind}_G(g))^3$ for some
$g \in G$. Furthermore, they put forward the following conjecture.\\[1ex]
{\bf Conjecture} \cite{MR2231893} Let $\chi$ be a primitive irreducible character of
an arbitrary finite group $G$. Then $\chi(1)$ divides
$\mbox{Ind}_G(g)$ for some element $g \in G$.\\[1ex]
They checked that the conjecture holds for all irreducible characters
(primitive or not) of all groups in the
Atlas \cite{MR827219}.

\subsection{$k(GV)$-problem}
In a slightly different direction Brauer \cite{MR0178056} was
interested in the number of characters in a given block. Let $p$ be
a prime and consider the number of  ordinary irreducible characters
belonging to the $p$-block $B$ with defect group $D$. He asked
whether this is less than or equal to $|D|$. For definitions of
blocks and related concepts see \cite[Chapter 15]{MR2270898}.

We know that the total number of irreducible characters is the number of conjugacy classes. Interestingly
in the $p$-soluble case the conjecture reduces to the following
problem \cite{MR0152569}: \\[1ex]
{\bf $k(GV)$-problem} Let $G$ be a finite $p'$-group for some prime
$p$ and let $V$ be a faithful ${\mathbb F}_p$-module for $G$.
Show that the number of conjugacy classes of the semidirect product $VG$ is bounded
by $|V |$.\\[1ex]
In this situation there is only one $p$-block so that all characters
are in the same $p$-block with defect group $V$ and Brauer's
question is about counting the number of conjugacy classes in a
group.

 The fundamental ideas for attacking the problem were developed around
1980 by Kn{\"o}rr~\cite{MR740615}, who proved the conjecture for $G$
supersoluble and later for $|G|$ odd \cite{MR748226}, the first
significant partial results. An important step on the way was the paper by Robinson \& Thompson, \cite{MR1407890},
which solved the problem for all but finitely many primes.

 The problem was finally solved by
Gluck, Magaard, Riese \& Schmid in 2004, \cite{MR2078936}. This was
the final piece of the jigsaw after a long series of papers. A more detailed
discussion and details of proofs are in the book by Schmid~\cite{MR2369608}. One of
the consequences of these results is that if $G$ is a $p'$-group
that can be embedded in GL$(m,p)$ for some integer $m$ then
$k(G)\leq p^m-1$.
\subsection{Huppert's Conjectures}
Let $G$ be a finite soluble group. Let $\bar{\sigma}(G)$ denote the
number of primes dividing any one character degree of $G$ and let
$\bar{\rho}(G)$ be the set of all primes which divide some character
degree of $G$. In 1985, Huppert~\cite{MR1112156} conjectured that if $G$ is soluble then $|\bar{\rho}(G)|
\leq 2 \bar{\sigma}(G)$. This conjecture is still open, the best result so
far, for soluble groups, is that $|\bar{\rho}(G)| \leq 3 \bar{\sigma}(G) +2$
\cite{MR898721}.

He also asked analogously whether $|\rho^*(G)| \leq 2 \sigma^*(G)$?
 Baer's results referred to earlier, Subsection~\ref{baerp}, show
that if $|\sigma^*(G)|=1$ then $|\rho^*(G)|\leq 2$. In \cite{MR1682749}
Mann considered groups in which all conjugacy class sizes involve at
most 2 primes, that is groups $G$ for which $|\sigma^*(G)|=2$. Such a
group is either soluble, or $G = Z(G) \times S$, where $S$ is isomorphic
to either $A_5$ or $SL(2,8)$ and $\rho^*(G)\leq 4$. This result
also appears in \cite{MR1256158} where the author is concentrating
on Huppert's Conjecture, note that both Baer's and Mann's papers
predate Huppert's conjecture. In \cite{MR1627921}  Zhang shows that
$\rho^*(G)\leq 4\sigma^*(G)$ for all soluble groups which
generalises the results of some earlier authors,
\cite{MR1201921,MR1128643,MR1152431,MR1049135}. Casolo in
\cite{MR1120130} proved the following theorem:
\begin{theo} Let $G$ be a  group which is $p$-nilpotent for at most one prime divisor of $|G|$ (this family includes all nonabelian simple groups), then
$$\rho^*(G)\leq 2\sigma^*(G).$$ \end{theo}
So Huppert's conjecture holds for many groups, consequently the following
theorem is somewhat surprising:
\begin{theo} \cite{MR1438293} There exist an infinite sequence of finite supersoluble, metabelian groups
$\{G_n\}$, such that $|\rho^*(G_n)|/\sigma^*(G_n)$ tends to 3 as $n
\to \infty$.\end{theo} So Huppert's conjecture is false. However,
they do show that $|\rho^*(G)| \leq 4 \sigma^*(G) + 2$ for all
soluble $G$.\label{sigma} Despite the result being false Casolo \& Dolfi have shown that there are linear bounds in both the character and conjugacy class versions of the conjecture for all groups.
\begin{theo}\cite{MR2352031} Let $G$ be a finite group. Then\\
(a) $|\bar{\rho}(G)|\leq 7\bar{\sigma}(G)$,\\
(b) $|\rho^*(G)|\leq 5\sigma^*(G).$
\end{theo}
In 2000 Huppert made the following conjecture:\\[1ex]
{\bf Conjecture} \cite{MR1804317} If $H$ is any simple nonabelian finite group and $G$ is a finite group such that ${\rm cd}(G)={\rm cd}(H)$, then $G\cong H\times A$, where $A$ is abelian.

The complete answer to this question is still to be found in the
character case. So we put forward the following question:-
\begin{ques} If $H$ is any simple nonabelian finite group and $G$ is a finite
group such that ${\rm cs}(G)={\rm cs}(H)$, then is $G\cong H\times A$,
where $A$ is abelian?\end{ques}

\newcommand{\etalchar}[1]{$^{#1}$}
\def\cprime{$'$} \def\cprime{$'$}

\end{document}